\documentclass[11pt]{article}
\usepackage{amssymb,amsmath,amsthm,fullpage}

 \newtheorem{thm}{Theorem}
 \newtheorem{cor}[thm]{Corollary}
 \newtheorem{lem}[thm]{Lemma}
 \newtheorem{prop}[thm]{Proposition}
 \newtheorem{conj}[thm]{Conjecture}
 \theoremstyle{definition}
 \newtheorem{defn}[thm]{Definition}
 \theoremstyle{remark}

\begin{document}

\title{A note on Reed's conjecture}

\author{ landon rabern \\
         rabernsorkin@cox.net }

\date{\today}

\maketitle

\begin{abstract}
In \cite{reed97}, Reed conjectures that the inequality $\chi (G) \leq \left \lceil \textstyle \frac{1}{2} (\omega (G) + \Delta (G) + 1) \right \rceil$ holds
for any graph $G$.  We prove this holds for a graph $G$ if $\overline{G}$ is disconnected.  From this it follows that the conjecture holds for
graphs with $\chi(G) > \left \lceil \frac{|G|}{2} \right \rceil$.  In addition, the conjecture holds for graphs with $\Delta(G) \geq |G| - \sqrt{|G| + 2\alpha(G) + 1}$.
In particular, Reed's conjecture holds for graphs with $\Delta(G) \geq |G| - \sqrt{|G| + 7}$.  Using these results, we proceed to show that if $|G|$ is an even order counterexample to Reed's conjecture, then $\overline{G}$ has a $1$-factor.  
Hence, for any even order graph
$G$, if $\chi(G) > \textstyle \frac{1}{2}(\omega(G) + \Delta(G) + 1) + 1$, then $\overline{G}$ is matching covered.
\end{abstract}

In all that follows, \emph{graph} will mean finite simple graph with non-empty vertex set.
Let $\mathbb{G}$ denote the collection of all graphs.  
Let $R_t \subseteq \mathbb{G}$ be the graphs satisfying $\chi \leq \textstyle \frac{1}{2}(\omega + \Delta + 1) + t$.
\begin{defn}
Given graphs $A$ and $B$, their \emph{join} $A+B$ is the graph with vertex set $V(A) \cup V(B)$ and edge set $E(A) \cup E(B) \cup \{ab \mid a \in V(A), b \in V(B) \}$.
Also, if $X$ and $Y$ are collections of graphs, we let $X+Y = \{A+B \mid A \in X, B \in Y \}$.
\end{defn}

First a few basic facts about joins.

\begin{lem}
Let $A$ and $B$ be graphs.  Then
\begin{itemize}
  \item[(a)] $|A + B| = |A| + |B|$,
  \item[(b)] $\omega(A+B) = \omega(A) + \omega(B)$,
  \item[(c)] $\chi(A+B) = \chi(A) + \chi(B)$,
  \item[(d)] $\Delta(A+B) = \max \{\Delta(A) + |B|, |A| + \Delta(B) \}$.
\end{itemize}
\end{lem}

\begin{proof}
These all follow immediately from the definitions.
\end{proof}

We will need the following result from \cite{rabern} and its immediate corollaries.

\begin{thm}
Let $I_1, \ldots, I_m$ be disjoint independent sets in a graph G. Then
\begin{equation}\label{indepdendentsetinequality} \chi (G) \leq \frac{1}{2} \left (\omega (G) + |G| - \displaystyle \sum_{j=1}^m |I_j| + 2m -1 \right). \end{equation}
\end{thm}

\begin{cor}
Let $G$ be a graph. Then
\[\chi (G) \leq \textstyle \frac{1}{2} (\omega (G) + |G| - \alpha(G) + 1).\]
\end{cor}
\begin{cor}
Let $G$ be a graph. Then
\[\chi (G) \leq \textstyle \frac{1}{2} (\omega (G) + |G|).\]
\end{cor}

\begin{prop}
$\mathbb{G} + R_t \subseteq R_t$ for all $t \in \mathbb{R}$.
\end{prop}
\begin{proof}
Fix $t \in \mathbb{R}$. Let $G \in \mathbb{G}$ and $H \in R_t$.  Applying Corollary 5 to $G$ gives

\[\chi (G) \leq \textstyle \frac{1}{2} (\omega (G) + |G|).\]

Also, since $H \in R_t$,

\[\chi(H) \leq \textstyle \frac{1}{2}(\omega(H) + \Delta(H) + 1) + t.\]

Adding these inequalities and applying Lemma 2 (b) and (c) gives

\[\chi(G+H) \leq \textstyle \frac{1}{2}(\omega(G+H) + |G| + \Delta(H) + 1) + t.\]
 
Now Lemma 2(d) gives $|G| + \Delta(H) \leq \Delta(G+H)$ and the result follows.
\end{proof}

\bigskip

The following two lemmas are special cases of Lemma 2 in \cite{molloy}.  
\begin{lem}
Let $G$ be a graph with $\chi(G) > \left \lceil \frac{|G|}{2} \right \rceil$.  Then there exists $X \subseteq V(G)$ such that
$\overline{G - X}$ is disconnected and $\chi(G - X) = \chi(G)$.
\end{lem}

\begin{lem}
If $G$ is a vertex critical graph with $\chi(G) > \left \lceil \frac{|G|}{2} \right \rceil$, then $\overline{G}$ is disconnected.
\end{lem}

\begin{prop}
If $G$ is a graph with $\alpha(G) \leq 2$, then $G \in R_{\frac{1}{2}}$.
\end{prop}
\begin{proof}
Assume this is not the case and let $G$ be a counterexample with the minimum number of vertices, say $|G| = n$.  
Since $\alpha(G) \leq 2$, we see that $V(G) \smallsetminus N(v) \cup \{v\}$ induces a clique for each $v \in V(G)$.
Hence $\omega(G) \geq n-\delta(G)-1$ which gives
\[ \Delta(G) + 1 \geq \delta(G) + 1 \geq n - \omega(G).\]

Since $G \not \in R_{\frac{1}{2}}$,

\begin{align*}
\chi(G) & > \textstyle \frac{1}{2} (\omega(G) + \Delta(G) + 1) + \frac{1}{2} \\
& \geq \textstyle \frac{1}{2} (\omega(G) + n - \omega(G)) + \frac{1}{2} \\
& = \frac{n+1}{2} \\
\end{align*}

Hence $\chi(G) > \left \lceil \frac{n}{2} \right \rceil$.  Now, using minimality of $G$, we see that $G$ is vertex critical.
Thus $\overline{G}$ is disconnected by Lemma 8. Hence we have $m \geq 2$ and non-empty graphs $C_1, \ldots, C_m$ such that $G = C_1 + \cdots + C_m$.
But, for $1 \leq i \leq m$, minimality of $G$ gives $C_i \in R_{\frac{1}{2}}$ since $\alpha(C_i) \leq \alpha(G) \leq 2$ and $|C_i| < n$.  Hence $G = C_1 + \cdots + C_m \in R_{\frac{1}{2}}$ 
by Proposition 6.  This contradiction completes the proof.
\end{proof}

\begin{defn}
Let $G$ be a graph and $r$ a positive integer.  A collection of disjoint independent sets in $G$ each with at least $r$ vertices
will be called an \emph{$r$-greedy} partial coloring of $G$.  A vertex of $G$ is said to be \emph{missed} by a partial coloring
just in case it appears in none of the independent sets.
\end{defn}

Applying Theorem 3 to an $r$-greedy partial coloring gives the following.
\begin{lem}
Let $G$ be a graph which is not complete and $M$ an $r$-greedy partial coloring of $G$.  Then
\[\chi(G) \leq \textstyle \frac{1}{2} (\omega(G) + |G| - (r-2)|M| - 1).\]
\end{lem}

\begin{lem}
Let $G$ be a graph and of all $3$-greedy partial colorings of $G$, let $M$ be one that misses the minimum number of vertices. Then
\begin{equation}\label{secondbound}
\chi(G) \leq \textstyle \frac{1}{2} (\omega(G) + \Delta(G) + 1) + \frac{|M|+1}{2}.
\end{equation}
\end{lem}
\begin{proof}
The first case to consider is when $M$ misses zero vertices.  In this case, $M$ is a proper coloring of $G$ and hence $\chi(G) \leq |M|$. Thus
\[ \chi(G) \leq \textstyle \frac{1}{2}(\chi(G) + |M|) \leq \frac{1}{2}(\Delta(G) + 1 + |M|) \leq \frac{1}{2} (\omega(G) + \Delta(G) + 1) + \frac{|M|+1}{2}.\]

Otherwise, $M$ misses at least one vertex and by the minimality condition placed on $M$, each vertex missed by $M$ must be adjacent to at least one vertex in each element of $M$.
Hence $\Delta(G - \cup M) \leq \Delta(G) - |M|$.  In addition, $\alpha(G - \cup M) \leq 2$.  Thus, applying Proposition 9 to $G - \cup M$, yields
\begin{align*}
\chi(G) &\leq |M| + \chi(G - \cup M) \\
&\leq |M| + \textstyle \frac{1}{2} (\omega(G - \cup M) + \Delta(G - \cup M) + 1) + \textstyle \frac{1}{2} \\
&\leq |M| + \textstyle \frac{1}{2} (\omega(G) + \Delta(G) - |M| + 1) + \textstyle \frac{1}{2} \\
&= \textstyle \frac{1}{2} (\omega(G) + \Delta(G) + 1) + \textstyle \frac{|M| + 1}{2}. \\
\end{align*}
\end{proof}

Taking $r = 3$ in Lemma 11 and adding the inequality with \eqref{secondbound} gives a better bound than Corollary 5.
\begin{prop}
Let $G$ be a graph.  Then $\chi(G) \leq \frac{1}{2} \left (\omega(G) + \frac{|G| + \Delta(G) + 1}{2} \right).$
\end{prop}

\begin{prop}
Let $A$ and $B$ be graphs. Then $A + B \in R_0$.
\end{prop}
\begin{proof}
Applying Proposition 13 to $A$ and $B$ and adding the inequalities yields

\[\chi(A) + \chi(B) \leq \frac{1}{2} \left (\omega(A) + \omega(B) + \frac{ \Delta(A) + |B| + |A| + \Delta(B) + 2}{2} \right ).\]

Using Lemma 2 (b),(c), and (d), this becomes

\[\chi(A + B) \leq \frac{1}{2} \left (\omega(A + B) + \frac{ 2 \Delta(A + B) + 2}{2} \right ) = \textstyle \frac{1}{2} (\omega(A + B) + \Delta(A + B) + 1).\]

Hence $A + B \in R_0$.
\end{proof}

\begin{cor}
If $G$ is a graph with $\chi(G) > \left \lceil \frac{|G|}{2} \right \rceil$, then $G \in R_0$.
\end{cor}
\begin{proof}
Let $G$ be a graph with $\chi(G) > \left \lceil \frac{|G|}{2} \right \rceil$.  Then, by Lemma 7, we have $X \subseteq V(G)$ such that
$\overline{G - X}$ is disconnected and $\chi(G - X) = \chi(G)$.  Since $\overline{G - X}$ is disconnected, there exist graphs $A$ and $B$ such that $G - X = A + B$.
Hence, by Proposition 14, 
\[\chi(G) = \chi(G - X) \leq \textstyle \frac{1}{2} (\omega(G - X) + \Delta(G - X) + 1) \leq \textstyle \frac{1}{2} (\omega(G) + \Delta(G) + 1).\]
Whence $G \in R_0$.

\end{proof}

\begin{cor}
Let $G$ be a graph and $t \geq 0$.  If $G \not \in R_t$, then $\Delta(G) + 1 \leq |G| - 2t - \omega(G)$.
\end{cor}
\begin{proof}
This is an immediate consequence of the previous corollary.
\end{proof}

\begin{lem}
Let $G$ be a graph with $\alpha(G) \leq 2$.  Then $\omega(G)^2 + \omega(G) \geq |G|$.
\end{lem}
\begin{proof}

Let $K$ be a maximal clique in $G$.  Then each vertex of $G - K$ is non-adjacent to at least one vertex in $K$ and hence some
vertex $v \in K$ is non-adjacent to at least $\frac{|G - K|}{|K|}$ vertices.  Since $\alpha(G) \leq 2$, the vertices non-adjacent to
$v$ form a clique.  Whence $\omega(G) \geq \frac{|G - K|}{|K|} = \frac{|G| - \omega(G)}{\omega(G)}$, which yields

\[\omega(G)^2 + \omega(G) \geq |G|. \]

\end{proof}

\begin{prop}
If $G$ is a graph with $\Delta(G) \geq |G| - \sqrt{|G| + 2\alpha(G) + 1}$, then $G \in R_{\frac{1}{2}}$.
\end{prop}
\begin{proof}
We prove the contrapositive.  Let $G$ be a graph with $n$ vertices, maximal degree $\Delta$, clique number $\omega$, and independence number $\alpha$ such
that $G \not \in R_{\frac{1}{2}}$.  Let $I$ be a maximal independent set in $G$.  Let $S$ be a maximal collection of disjoint $3$-vertex indepedent
sets of $G - I$.  Since $\alpha(G - (\cup S) \cup I) \leq 2$ , we may apply Lemma 17 to get $\omega^2 + \omega \geq |G - (\cup S) \cup I| = n - \alpha - 3|S|$. Hence 

\begin{equation}\label{SIsBig}
|S| \geq \frac{n - \alpha - (\omega^2 + \omega)}{3}.
\end{equation}
Now, using the fact that $G \not \in R_{\frac{1}{2}}$ with Lemma 11, we have $n - \alpha - |S| + 1 > \Delta + 2$.  Putting this together with
\eqref{SIsBig} we have

\[n - \alpha - \Delta - 1 > |S| \geq \frac{n - \alpha - (\omega^2 + \omega)}{3}.\]

Which implies that

\begin{equation}\label{DeltaBound}
\Delta < \frac{2n + \omega^2 + \omega - 2\alpha - 3}{3}.
\end{equation}

By Corollary 16, $\omega \leq n - \Delta - 2$.  Plugging this into \eqref{DeltaBound} and doing a little algebra, we find that
$\Delta < n - \sqrt{n + 2\alpha + 1}$.  This completes the proof of the contrapositive.
\end{proof}

\begin{cor}
If $G$ is a graph with $\Delta(G) \geq |G| - \sqrt{|G| + 7}$, then $G \in R_{\frac{1}{2}}$.
\end{cor}
\begin{proof}
Let $G$ be a graph with $\Delta(G) \geq |G| - \sqrt{|G| + 7}$.  If $\alpha(G) \leq 2$, then $G \in R_{\frac{1}{2}}$ by Proposition 9.  Otherwise,
$\alpha(G) \geq 3$ and $G \in R_{\frac{1}{2}}$ by Proposition 18.
\end{proof}

\begin{lem}
Let $G$ be a graph and $t \in \frac{1}{2} \mathbb{Z}$.  If $G \not \in R_t$, then $\Delta(G) + 1 \leq |G| - 2t - \alpha(G)$.
\end{lem}
\begin{proof}
Assume $G \not \in R_t$.  Applying Corollary 4 gives 
\[\textstyle \frac{1}{2} ( \omega(G) + \Delta(G) + 1 ) + t < \chi(G) \leq \textstyle \frac{1}{2} ( \omega(G) + |G| - \alpha(G) + 1 ).\]
The lemma follows.
\end{proof}

Note that for $t \geq \frac{1}{2}$ in the lemma, we must have $\alpha(G) \geq 3$ by Propostion 9, so the lemma gives $\Delta(G) + 1 \leq |G| - 2t - 3$.

\begin{lem}
If $k \geq 2$ and $G_1, \ldots, G_k$ are graphs with $\Delta(G_i) + 1 \leq |G_i| - 3$ for each $i$, then 
\[G_1 + \cdots + G_k \in R_{2-k}.\]
\end{lem}
\begin{proof}
Assume this is not the case and let $G_1, \ldots, G_k$ constitute a counterexample with the smallest $k$.  Then, by Proposition 14, $k > 2$.  
Set $D = G_1 + \cdots + G_{k - 1}$.  Note that $D \in R_{2-(k-1)} = R_{3-k}$ by the minimality of $k$.  
Let $t \in \frac{1}{2} \mathbb{Z}$ be minimal such that $G_k \in R_t$.  Since $G_k \not \in R_{t - \frac{1}{2}}$, using Lemma 20 for $t \geq 1$ and the
fact that $\Delta(G_k) + 1 \leq |G_k| - 3$ for $t \leq \frac{1}{2}$, we find that $\Delta(G_k) + 1 \leq |G_k| - 2t - 2$.
We have,
\begin{align*}
\chi(D + G_k) &= \chi(D) + \chi(G_k) \\
&\leq \textstyle \frac{1}{2} ( \omega(D) + \omega(G_k) + \Delta(D) + \Delta(G_k) + 2) + 3-k + t \\
&= \textstyle \frac{1}{2} ( \omega(D + G_k) + \Delta(D) + \Delta(G_k) + 2) + 3-k + t \\
&= \textstyle \frac{1}{2} ( \omega(D + G_k) + \Delta(D) + |G_k| + 1 + \Delta(G_k) - |G_k| + 1) + 3-k + t \\
&\leq \textstyle \frac{1}{2} ( \omega(D + G_k) + \Delta(D + G_k) + 1) + \frac{1}{2}(\Delta(G_k) - |G_k| + 1) + 3 - k + t \\
&\leq \textstyle \frac{1}{2} ( \omega(D + G_k) + \Delta(D + G_k) + 1) + \frac{1}{2}(-2t - 3 + 1) + 3 - k + t \\
&= \textstyle \frac{1}{2} ( \omega(D + G_k) + \Delta(D + G_k) + 1) + 2 - k. \\
\end{align*}

Hence $G_1 + \cdots + G_k \in R_{2-k}$, contradicting our assumption.
\end{proof}

The hypotheses of this lemma can be weakened, but we do not use the following stronger lemma in what follows.

\begin{lem}
If $k \geq 2$ and $G_1, \ldots, G_k$ are graphs, which are not $5$-cycles, with $\Delta(G_i) + 1 \leq |G_i| - 2$ for each $i$, then 
\[G_1 + \cdots + G_k \in R_{2-k}.\]
\end{lem}
\begin{proof}
Similar to Lemma 21.  Graphs with $\Delta(G_i) + 1 \leq |G_i| - 3$ only matter if $G_i \in R_{\frac{1}{2}} \smallsetminus R_{0}$.  Corollary 16 shows that
such graphs have $\omega(G_i) \leq 2$ and Corollary 4 shows they have $\alpha(G) \leq 2$.  Thus they have order less than $6$ and we see that the only one that 
breaks the lemma is $C_5$.
\end{proof}

\begin{defn}
The \emph{matching number} of a graph $G$, denoted $\nu (G)$ is the number of edges in a maximal matching of $G$.
\end{defn}

\begin{prop}
Let $G$ be a graph.  If $\nu(\overline{G}) < \left \lfloor \frac{|G|}{2} \right \rfloor$, then $G \in R_0$.
\end{prop}
\begin{proof}
Assume $\nu(\overline{G}) < \left \lfloor \frac{|G|}{2} \right \rfloor$.  Then, by Tutte's Theorem, we have $X \subseteq V(G)$ such that
$\overline{G - X}$ has at least $m$ odd components; where $m \geq |X| + 2$ if $|G|$ is even and $m \geq |X| + 3$ if $|G|$ is odd.  Hence, we have graphs $G_1, \ldots, G_{m}$ such that $G - X  = G_1 + \cdots + G_{m}$.
Note that by picking one vertex from each component we induce a clique.  Hence $\omega(G) \geq m$.
To get a contradiction, assume $G \not \in R_0$.
First assume there is some $G_i$ for which $\Delta(G_i) + 1 \geq |G_i| - 2$, then 
\[\Delta(G) + 1 \geq \Delta(G - X) + 1 \geq |G_1| + \cdots + |G_{i-1}| + \Delta(G_i) + |G_{i+1}| + \cdots + |G_{m}| + 1 \geq |G| - |X| - 2.\]
Since $G \not \in R_0$,
\[ \chi(G) > \textstyle \frac{1}{2} ( \omega(G) + \Delta(G) + 1) \geq  \frac{1}{2} ( m + (|G| - |X| - 2)) \geq \left \lceil \frac{|G|}{2} \right \rceil. \]
Hence $G \in R_0$ by Corollary 15!  Thus, we may assume $\Delta(G_i) + 1 \leq |G_i| - 3$ for each $i$.  Now Lemma 21 yields $G - X \in R_{-|X|}$.  
Whence $G \in R_0$.  This contradiction completes the proof.
\end{proof}

\begin{cor}
Let $G$ be an even order graph.  If $G \not \in R_0$, then $\overline{G}$ has a $1$-factor.
\end{cor}

\begin{defn}
A graph is called \emph{matching covered} if every edge participates in a perfect matching.
\end{defn}

\begin{cor}
Let $G$ be an even order graph with $G \not \in R_1$.  Then $\overline{G}$ is matching covered.
\end{cor}

Lemma 21 can be generalized.

\begin{prop}
Let $m \in \mathbb{N}$.  Let $k \geq 2$ and $G_1, \ldots, G_k$ be graphs such that $|G_i| < r(m,m) \Rightarrow G_i \in R_{\frac{1}{2}}$. 
If $\Delta(G_i) + 1 \leq |G_i| - m$ for each $i$, then \[G_1 + \cdots + G_k \in R_{(m-1)(1 - \frac{k}{2})}.\]
\end{prop}
\begin{proof}
Assume this is not the case and let $G_1, \ldots, G_k$ constitute a counterexample with the smallest $k$.  Then, by Proposition 14, $k > 2$.  
Set $D = G_1 + \cdots + G_{k - 1}$.  Note that $D \in R_{(m-1)(1 - \frac{k-1}{2})}$ by the minimality of $k$.  
Let $t \in \frac{1}{2} \mathbb{Z}$ be minimal such that $G_k \in R_t$.  We would like to have $\Delta(G_k) + 1 \leq |G_k| - 2t - (m-1)$.
If $t \leq \frac{1}{2}$, then we are all good since $\Delta(G_k) + 1 \leq |G_k| - m$.  So, to get a contradiction, assume $t \geq 1$ and 
$\Delta(G_k) + 1 > |G_k| - 2t - (m-1)$.  Then, by Corollary 16, $\omega(G_k) \leq (m-1)$. Also, by Lemma 20, $\alpha(G_k) \leq (m-1)$. 
Hence $|G_k| < r(m,m)$ contradicting the fact that $G_k \not \in R_{\frac{1}{2}}$.  Hence we do indeed have $\Delta(G_k) + 1 \leq |G_k| - 2t - (m-1)$.

We have,
\begin{align*}
\chi(D + G_k) &= \chi(D) + \chi(G_k) \\
&\leq \textstyle \frac{1}{2} ( \omega(D) + \omega(G_k) + \Delta(D) + \Delta(G_k) + 2) + (m-1)(1 - \frac{k-1}{2}) + t \\
&= \textstyle \frac{1}{2} ( \omega(D + G_k) + \Delta(D) + \Delta(G_k) + 2) + (m-1)(1 - \frac{k-1}{2}) + t \\
&= \textstyle \frac{1}{2} ( \omega(D + G_k) + \Delta(D) + |G_k| + 1 + \Delta(G_k) - |G_k| + 1) + (m-1)(1 - \frac{k-1}{2}) + t \\
&\leq \textstyle \frac{1}{2} ( \omega(D + G_k) + \Delta(D + G_k) + 1) + \frac{1}{2}(\Delta(G_k) - |G_k| + 1) + (m-1)(1 - \frac{k-1}{2}) + t \\
&\leq \textstyle \frac{1}{2} ( \omega(D + G_k) + \Delta(D + G_k) + 1) + \frac{1}{2}(-2t - 4 + 1) + (m-1)(1 - \frac{k-1}{2}) + t \\
&= \textstyle \frac{1}{2} ( \omega(D + G_k) + \Delta(D + G_k) + 1) + (m-1)(1 - \frac{k}{2}). \\
\end{align*}

Hence $G_1 + \cdots + G_k \in R_{(m-1)(1 - \frac{k}{2})}$, contradicting our assumption.
\end{proof}

\begin{conj}
Let $m \in \mathbb{N}$.  If $k \geq 2$ and $G_1, \ldots, G_k$ are graphs with $\Delta(G_i) + 1 \leq |G_i| - m$ for each $i$, 
then \[G_1 + \cdots + G_k \in R_{(m-1)(1 - \frac{k}{2})}.\]
\end{conj}

We can do a bit better than Proposition 28 in the following special case.

\begin{lem}
If $k \geq 2$ and $G_1, \ldots, G_k$ are non-complete graphs, then 
\[G_1 + \cdots + G_k \in R_{1-\frac{k}{2}}.\]
\end{lem}
\begin{proof}
Assume this is not the case and let $G_1, \ldots, G_k$ constitute a counterexample with the smallest $k$.  Then, by Proposition 14, $k > 2$.  
Set $D = G_1 + \cdots + G_{k - 1}$.  Note that $D \in R_{1-\frac{(k-1)}{2}}$ by the minimality of $k$.  
Let $t \in \frac{1}{2} \mathbb{Z}$ be minimal such that $G_k \in R_t$.  Since $G_k \not \in R_{t - \frac{1}{2}}$, if $t \geq \frac{1}{2}$, then, by Corollary 16,
$\Delta(G_k) + 1 \leq |G_k| - 2t - 1$.  If $t \leq 0$ and $\Delta(G_k) + 1 > |G_k| - 2t - 1$, then $t = 0$ and $\Delta(G_k) + 1 = |G_k|$; however, Corollary 4 shows
that the only such graphs are complete graphs which we have excluded.  Whence $\Delta(G_k) + 1 \leq |G_k| - 2t - 1$.
We have,
\begin{align*}
\chi(D + G_k) &= \chi(D) + \chi(G_k) \\
&\leq \textstyle \frac{1}{2} ( \omega(D) + \omega(G_k) + \Delta(D) + \Delta(G_k) + 2) + 1-\frac{(k-1)}{2} + t \\
&= \textstyle \frac{1}{2} ( \omega(D + G_k) + \Delta(D) + \Delta(G_k) + 2) + 1-\frac{(k-1)}{2} + t \\
&= \textstyle \frac{1}{2} ( \omega(D + G_k) + \Delta(D) + |G_k| + 1 + \Delta(G_k) - |G_k| + 1) + 1-\frac{(k-1)}{2} + t \\
&\leq \textstyle \frac{1}{2} ( \omega(D + G_k) + \Delta(D + G_k) + 1) + \frac{1}{2}(\Delta(G_k) - |G_k| + 1) + 1-\frac{(k-1)}{2} + t \\
&\leq \textstyle \frac{1}{2} ( \omega(D + G_k) + \Delta(D + G_k) + 1) + \frac{1}{2}(-2t - 2 + 1) + 1-\frac{(k-1)}{2} + t \\
&= \textstyle \frac{1}{2} ( \omega(D + G_k) + \Delta(D + G_k) + 1) + 1-\frac{k}{2}. \\
\end{align*}

Hence $G_1 + \cdots + G_k \in R_{1-\frac{k}{2}}$, contradicting our assumption.
\end{proof}

Similar ideas can be used to prove theorems in the same vein as Proposition 24 with $R_0$ replaced by $R_{-a}$ for $a > 0$.  However, the details get
hairy and we don't feel they are worth reproducing here as they don't seem to give new insight into Reed's conjecture.

\end{document}